\documentclass{notices}

\usepackage{amsmath,amssymb,amsfonts,amsthm}
\usepackage{mathtools}
\usepackage{graphicx}
\usepackage{tabularx}
\usepackage{array}
\usepackage{booktabs}
\usepackage{xcolor}
\usepackage{tikz}
\usetikzlibrary{arrows.meta,calc,positioning,decorations.pathreplacing,patterns,fit,backgrounds}

\definecolor{kinblue}{RGB}{33,91,126}
\definecolor{kinred}{RGB}{156,58,52}
\definecolor{kingreen}{RGB}{72,119,80}
\definecolor{kinateal}{RGB}{42,128,125}
\definecolor{kinorange}{RGB}{191,117,44}
\definecolor{kinpaper}{RGB}{247,245,238}
\definecolor{kinink}{RGB}{45,45,45}

\tikzset{
  noticebox/.style={draw=kinink!55, rounded corners=6pt, line width=0.45pt,
    fill=kinpaper, inner sep=6pt, align=center},
  noticearrow/.style={-Latex, line width=0.8pt, draw=kinink!80},
  noticelabel/.style={font=\footnotesize, text=kinink!90, align=center}
}

\DeclareRobustCommand{\rev}[1]{#1}
\newcommand{\revm}[1]{#1}

\newtheorem{principle}{Principle}[section]
\newtheorem{theorem}[principle]{Theorem}

\newenvironment{ideaproof}{\begin{proof}[Idea of the proof]}{\end{proof}}

\newcommand{\eps}{\varepsilon}
\newcommand{\T}{\mathbb T}
\newcommand{\R}{\mathbb R}

\newcommand{\dd}{\,\mathrm d}

\title{Quasineutral Plasmas and the Geometry of Kinetic Stability}

\date{}

\author{
Mikaela Iacobelli
\affil{
Mikaela Iacobelli is professor of mathematics at ETH Z\"urich.
Her email address is \textsf{mikaela.iacobelli@math.ethz.ch}.
}
}

\begin{document}

\maketitle

\section*{A kinetic view of plasmas}

If we heat up a neutral gas or subject it to a strong electromagnetic field,
some electrons may detach from the atoms to which they were bound.  This
process is called ionization.  The gas then contains free electrons and
positively charged ions, and it may also contain particles that have remained
neutral.  When the charged component is large enough to influence the
\rev{dynamics one wants to describe}, one calls the medium a plasma.

Plasmas are found in stars, in the solar wind, and in the interstellar medium.
They can also be produced in the laboratory, for instance, in devices designed
for controlled fusion.  The degree of ionization may vary considerably from one
example to another; thus, there is no universal threshold at which an ionized
gas becomes a plasma, and the relevant criterion is whether the charged
particles play a leading role in the dynamics.

This is the feature that distinguishes plasmas from neutral gases.  In a
neutral dilute gas, the interaction between particles is often modeled through
collisions: particles travel freely until they meet and exchange
momentum.  In a plasma, charged particles interact predominantly through the
electromagnetic field generated by their distribution.  A local imbalance of
charge creates a field, the field acts on the particles, and the motion of the
particles changes the distribution of charge again.  This self-consistent
interaction is at the heart of the nonlinear mechanism that governs the
evolution of a plasma, and it is the starting point for the kinetic models
discussed below.

The point of view adopted in this article is that of kinetic theory, a branch
of statistical mechanics developed to describe systems made of many particles.
Instead of following the trajectory of each particle separately, one studies the evolution of a distribution function, which encodes the state
of the plasma at a statistical level. 
This function depends on time, position, and velocity, and it describes how the particles are distributed in space, and
also how they are distributed among different velocities.

The unknown is then a distribution function
\[
   f(t,x,v), \qquad (x,v)\in \T^d\times\R^d,
\]
on the one-particle phase space.  The quantity \(f(t,x,v)\,dx\,dv\) represents
the fraction of particles which, at time \(t\), are near the position \(x\) and
have velocity near \(v\).

This description is intermediate between particles and fluids.  At the
microscopic level, one would follow all individual particles, with their
positions and velocities, but this task becomes impossible when the number of
particles \rev{is} of the order of Avogadro's number, and in any case, it would give
far more detailed information than one can use.  At the fluid level, one looks
instead at a much coarser picture, following the evolution in time of
macroscopic observables, which are averaged quantities such as the density and
the mean velocity.  The kinetic description lies between these two points of
view: it does not follow each particle separately, but it still keeps the
velocity variable, and therefore retains information that is lost in a macroscopic, 
fluid model.

\begin{figure*}[t]
\centering
\resizebox{0.96\textwidth}{!}{%
\begin{tikzpicture}[>=Latex, font=\small]
\node[noticebox, minimum width=3.55cm, minimum height=2.95cm] (micro) at (0,0) {};
\node[noticebox, minimum width=3.95cm, minimum height=2.95cm] (kinetic) at (5.25,0) {};
\node[noticebox, minimum width=3.95cm, minimum height=2.95cm] (macro) at (10.65,0) {};

\node[font=\bfseries\small] at (0,1.05) {microscopic};
\node[noticelabel] at (0,0.63) {$N$ trajectories};
\foreach \x/\y/\c in {-1.10/-0.50/kinblue,-0.68/-0.10/kinred,-0.18/-0.32/kingreen,0.36/0.04/kinorange,0.86/-0.42/kinateal}
  {
  \fill[\c!70] (\x,\y) circle (0.075);
  \draw[-Latex, \c!80, line width=0.55pt] (\x,\y) to ++(0.24,0.16);
  }
\node[noticelabel] at (0,-1.12) {Newton equations};

\node[font=\bfseries\small] at (5.25,1.05) {kinetic};
\node[noticelabel] at (5.25,0.63) {distribution $f(t,x,v)$};
\draw[->, kinink!70] (4.30,-0.88) to (6.25,-0.88) node[right, font=\scriptsize] {$x$};
\draw[->, kinink!70] (4.52,-1.08) to (4.52,0.22) node[above, font=\scriptsize] {$v$};
\draw[fill=kinblue!12, draw=kinblue!85, line width=0.75pt]
  (4.75,-0.67) .. controls (4.90,-0.15) and (5.48,0.13) .. (5.92,-0.13)
  .. controls (6.20,-0.37) and (5.88,-0.75) .. (5.25,-0.73)
  .. controls (4.96,-0.72) and (4.80,-0.70) .. (4.75,-0.67);
\node[noticelabel] at (5.25,-1.28) {one-particle phase space};

\node[font=\bfseries\small] at (10.65,1.05) {macroscopic};
\node[noticelabel] at (10.65,0.63) {averaged quantities};
\node[font=\normalsize] at (10.65,-0.23) {$\rho,\ j,\ E$};
\draw[kinink!70, dashed, line width=0.7pt]
  (9.53,-0.70) .. controls (9.95,-0.38) and (10.30,-0.96) .. (10.65,-0.70)
  .. controls (11.06,-0.37) and (11.42,-0.95) .. (11.79,-0.70);
\node[noticelabel] at (10.65,-1.22) {coarser variables};

\draw[noticearrow] (1.90,0.22) to node[above, noticelabel] {$N\to\infty$}
  node[below, noticelabel] {mean field} (3.15,0.22);
\draw[noticearrow] (7.35,0.22) to node[above, noticelabel] {$\lambda_D/L\to0$}
  node[below, noticelabel] {QN limit} (8.55,0.22);
\end{tikzpicture}}
\caption{Three levels of description.  The mean-field limit replaces many
particles by a distribution \(f(t,x,v)\).  The quasineutral limit discussed in
this article starts from the kinetic equation when the ratio between the Debye
length and the observation scale satisfies \(\lambda_D/L\to0\); it leads to a
macroscopic constraint on the density, while still retaining the velocity
distribution.}
\label{fig:three-scales}
\end{figure*}
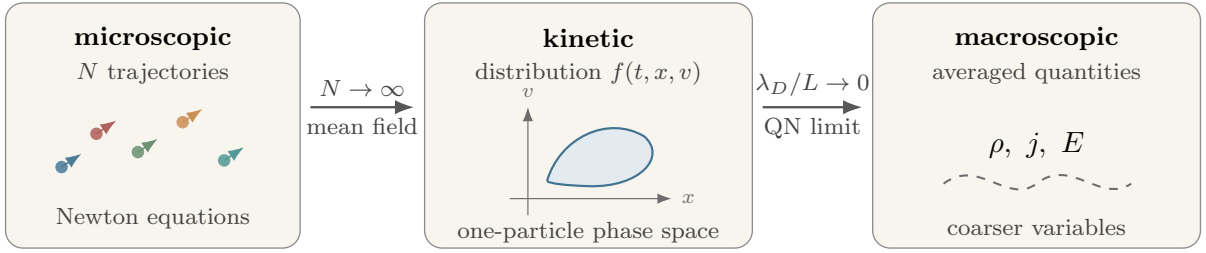

Let us first look at the most classical kinetic model for plasmas, the
Vlasov-Poisson (VP) system.  In the simplest electrostatic setting, one follows
one charged species, say the electrons, while the ions are modeled as a fixed
and uniformly distributed background.  \rev{Throughout the article, unless
otherwise specified, we take the position variable to be periodic,
\(x\in\T^d=\R^d/\mathbb Z^d\), in order not to discuss boundary effects.}  After a
normalization of the total mass
and of the background charge, the system on the torus reads
\[
\left\{
\begin{aligned}
   &\partial_t f+v\cdot\nabla_x f-\nabla_x U\cdot\nabla_v f=0,\\
   &-\Delta_x U=\rho_f-1,\\
   &\rho_f(t,x)=\int_{\R^d} f(t,x,v)\,dv.
\end{aligned}
\right.
\]
Here \(x\in\T^d\) is the position variable and \(v\in\R^d\) is the velocity
variable.   The function \(U(t,x)\) is the electrostatic potential, while
\[
   E(t,x)=-\nabla_x U(t,x)
\]
is the electric field.  The constant \(1\) in the Poisson equation represents
the neutralizing ionic background, so that the field is generated by the
deviation of the electron density \(\rho_f\) from this background.

The first equation says that the distribution function is transported in phase
space by the electric field.  Equivalently, the particles follow the
characteristics
\[
   \dot X(t)=V(t),\qquad
   \dot V(t)=-\nabla_x U(t,X(t)).
\]
The second equation closes the system by determining the field from the
density of particles and the self-consistent nature of the model.

But how did we arrive at this kinetic equation?  If we go back to the
microscopic description and consider the evolution of \(N\) charged particles
through Newton's laws of motion, we have to follow the position and velocity of
each particle.  In a simple mean-field model, this gives a system of the form
\[
   \dot x_i=v_i,\qquad
   \dot v_i=\frac1N\sum_{j\ne i}K(x_i-x_j),
   \qquad i=1,\ldots,N,
\]
where \(K\) is the interaction force.  \rev{In the electrostatic plasma model,
this force is the Coulomb force, up to the normalization chosen for the
background.} To this configuration, one can associate
the empirical measure
\[
   \mu_t^N=\frac1N\sum_{i=1}^N\delta_{(x_i(t),v_i(t))}.
\]
The idea of the mean-field limit is that, when the number of particles becomes
very large, this empirical measure should behave like a smooth distribution
function \(f(t,x,v)\), and the particle system should be replaced by a kinetic
equation.  This passage from particles to a density is a deep problem in its
own right, especially for singular interactions such as the Coulomb force; see,
for instance, \cite{DuerinckxSerfaty}.

Rigorously establishing the links between the microscopic, kinetic, and
fluid points of view is a central question in kinetic theory.  \rev{It gives
mathematical meaning to the idea that the various models are consistent with
each other and can be used at different scales, depending on the question.}

In this article, we start from the kinetic description and focus on a second
passage: the derivation of a fluid-like plasma model from the
(VP) system.  

In the following sections we will focus on a physical length scale that plays a
central role in plasma physics: the Debye length.  The size of this length,
compared with the scale at which the plasma is observed, determines whether the
full kinetic description is needed, or whether a more constrained, fluid-like
behavior can emerge.  The limiting procedure in which this ratio becomes small
is called the quasineutral limit.  We will see that this limit is closely tied
to the stability and instability mechanisms of the plasma, and that one needs
quantitative ways of measuring how perturbations of the distribution function
are propagated by the dynamics.  The main mathematical tool discussed in this
article will be a class of anisotropic, or kinetic, Wasserstein distances.

%%%%%%%%%%%%%%

\section{When the Debye length is small}
We now start from the (VP) system and focus on the passage from a
kinetic description to a more fluid-like one; the physical scale governing this
passage is the Debye length, the characteristic length scale for charge
separation in a plasma.

Roughly speaking, if a small region carries an excess of charge, the
mobile charges in the plasma rapidly rearrange themselves and screen this defect over a
distance of the order of the Debye length, denoted by \(\lambda_D\).  For an
electron plasma one may write
\[
   \lambda_D=
   \left(\frac{\varepsilon_0 k_B T_e}{n_e e^2}\right)^{1/2},
\]
where \(T_e\) is the electron temperature and \(n_e\) is the electron density.
For us, the precise expression is less important than the fact that this
length depends on the physical state of the plasma, through quantities such as
temperature and density; a standard reference for the physical discussion of
screening is Chen's book \cite{Chen}.

The relevant quantity is then the ratio between the Debye length and the length
scale at which the plasma is observed.  If this observation scale is denoted by
\(L\), we introduce the dimensionless parameter
\[
   \eps=\frac{\lambda_D}{L}.
\]
The quasineutral regime corresponds to \(\eps\ll1\), meaning that the plasma is
observed on a scale much larger than the scale on which charge imbalances are
screened.  At this macroscopic scale the plasma appears almost neutral.

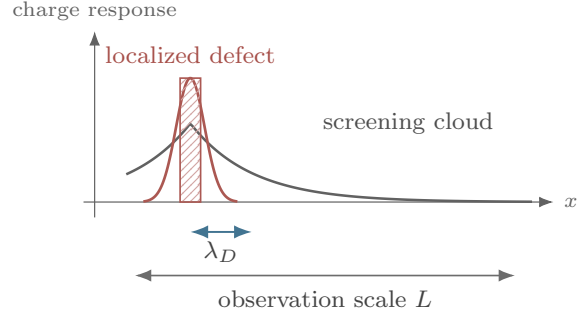
\begin{figure}[t]
\centering
\resizebox{0.96\columnwidth}{!}{%
\begin{tikzpicture}[>=Latex, font=\small]
\draw[->, kinink!75] (-0.25,0) to (5.25,0) node[right, font=\scriptsize] {$x$};
\draw[->, kinink!75] (-0.12,-0.18) to (-0.12,2.0) node[above, font=\scriptsize] {charge response};
\draw[kinink!75, line width=0.85pt, domain=0.25:5.0, samples=160]
  plot (\x,{0.92*exp(-abs(\x-1.0)/0.72)});
\draw[kinred!85, line width=0.9pt, domain=0.45:1.55, samples=90]
  plot (\x,{1.45*exp(-18*(\x-1.0)*(\x-1.0))});
\fill[kinred!15, pattern=north east lines, pattern color=kinred!45]
  (0.88,0) rectangle (1.12,1.45);
\draw[kinred!85, line width=0.55pt] (0.88,0) rectangle (1.12,1.45);
\node[noticelabel, text=kinred!90] at (1.0,1.72) {\rev{localized defect}};
\node[noticelabel, text=kinink!85] at (3.55,0.92) {screening cloud};
\draw[<->, kinblue!85, line width=0.65pt] (1.0,-0.35) to (1.72,-0.35);
\node[noticelabel] at (1.36,-0.58) {$\lambda_D$};
\draw[<->, kinink!70, line width=0.55pt] (0.35,-0.86) to (4.8,-0.86);
\node[noticelabel] at (2.58,-1.13) {observation scale $L$};
\end{tikzpicture}}
\caption{\rev{A localized charge defect, represented by the shaded red region,
is screened on the Debye scale \(\lambda_D\); the decaying profile sketches the
charge response around it.}  The quasineutral regime corresponds to
observations at a much larger scale \(L\).}
\label{fig:debye-screening}
\end{figure}

Let us now insert this scaling into the (VP) equation.  After
nondimensionalization, the electron (VP) system with immobile ions
takes the form
\[
(VP)_\eps
\qquad
\left\{
\begin{aligned}
\partial_t f_\eps+&v\cdot\nabla_x f_\eps
    +E_\eps\cdot\nabla_v f_\eps=0,\\
E_\eps&=-\nabla_x U_\eps,\\
-\eps^2\Delta_xU_\eps&=\rho_{f_\eps}-1,\\
\rho_{f_\eps}(t,x)&=\int_{\R^d}f_\eps(t,x,v)\,dv .
\end{aligned}
\right.
\tag{1}
\]
Here all variables are dimensionless, and the small parameter
\(\eps=\lambda_D/L\) is the ratio between the Debye length and the observation
scale; the quasineutral limit corresponds to \(\eps\to0\).

For every fixed \(\eps>0\), the last equation is an elliptic equation.  After
choosing a normalization of the potential, it allows one to recover
\(U_\eps\) from the charge imbalance \(\rho_{f_\eps}-1\) by inverting the
Laplacian.  Thus the (VP) system contains an elliptic coupling
between the density and the electric field; this coupling gives spatial
regularity to the field, at the level of the inverse Laplacian, and it is one
of the structural features of the original system.

As \(\eps\) tends to zero, this elliptic relation degenerates.  Written in the
form
\[
   \rho_{f_\eps}-1=-\eps^2\Delta_x U_\eps,
\]
the Poisson equation shows that, as long as the potential does not develop
second derivatives of size \(\eps^{-2}\), the charge imbalance must disappear
at the scale of observation.  The formal leading-order relation is therefore
\[
   \rho_f=1.
\]
This is the quasineutral constraint.  It says that, in the limit, the electron
density balances the background charge; at the same time, the potential is no
longer obtained from a Poisson equation driven by the charge imbalance.  It has
to be selected in a different way, so that the kinetic evolution preserves the
constraint \(\rho_f=1\).  This loss of the original elliptic coupling is one of
the reasons why the limiting system is more singular than the scaled
(VP) equations.

The formal limiting system is then
\[
(KInE)
\qquad
\left\{
\begin{aligned}
\partial_t f&+v\cdot\nabla_x f+E\cdot\nabla_v f=0,\\
E&=-\nabla_x U,\\
\rho_f&=1.
\end{aligned}
\right.
\tag{2}
\]
This system is called the kinetic incompressible Euler equation.  The
terminology reflects the analogy with incompressible fluid mechanics.  In the
Euler equation, the pressure is chosen so that the constraint
\(\nabla_x\cdot u=0\) is preserved; here, the electric potential \(U\) is
chosen so that the kinetic evolution preserves \(\rho_f=1\).  The analogy is
only partial, since the unknown remains the full phase-space distribution
\(f(t,x,v)\), rather than a single velocity field.

The way in which the potential is chosen in the limit can be seen at the level
of moments.  Define
\[
\begin{aligned}
   \rho_\eps&=\int_{\R^d} f_\eps\,dv,\qquad
   j_\eps=\int_{\R^d} v f_\eps\,dv,\\
   \Pi_\eps&=\int_{\R^d} v\otimes v\,f_\eps\,dv .
\end{aligned}
\]
Thus \(\rho_\eps\) is the density, \(j_\eps\) is the current, and \(\Pi_\eps\)
is the second velocity moment.  Integrating the Vlasov equation with respect
to \(v\), first with weight \(1\) and then with weight \(v\), gives formally
\[
   \partial_t\rho_\eps+\nabla_x\cdot j_\eps=0,
   \qquad
   \partial_t j_\eps+\nabla_x\cdot\Pi_\eps=\rho_\eps E_\eps.
\tag{3}
\]
In the limit \(\eps\to0\), the constraint \(\rho_f=1\) implies, through the
continuity equation,
\[
   \nabla_x\cdot j_f=0.
\]
The limiting momentum equation reads
\[
   \partial_t j_f+\nabla_x\cdot\Pi_f=-\nabla_x U,
   \qquad
   \Pi_f=\int_{\R^d}v\otimes v\,f\,dv .
\]
Taking its divergence and using \(\nabla_x\cdot j_f=0\), one obtains
\[
   -\Delta_x U=\nabla_x\cdot\nabla_x\cdot\Pi_f.
\tag{4}
\]
This is no longer the original Poisson equation, where the potential is
obtained from the charge imbalance.  It is instead a compatibility condition:
the potential is selected so that the constraint \(\rho_f=1\) remains true
along the kinetic evolution.

To recover a fluid equation, one has to add a closure.  A particularly simple
case is the monokinetic ansatz
\[
   f(t,x,v)=\rho(t,x)\,\delta_{v=u(t,x)},
\tag{5}
\]
which means that, at each point in space, all particles have the same velocity.
Then
\[
   j=\rho u,\qquad \Pi=\rho u\otimes u.
\]
Since the quasineutral constraint gives \(\rho=1\), the continuity equation
becomes
\[
   \nabla_x\cdot u=0,
\]
and the momentum equation reduces to
\[
   \partial_t u+u\cdot\nabla_x u=-\nabla_x U.
\tag{6}
\]
This is the incompressible Euler equation, with \(U\) playing the role of the
pressure \cite{Masmoudi}\rev{; see also \cite{Brenier2000}}.  Thus the Euler system appears as a closed fluid form of the
quasineutral limit when the velocity distribution collapses to one velocity at
each spatial point.  If several velocities coexist at the same point, the
limit retains its kinetic character.

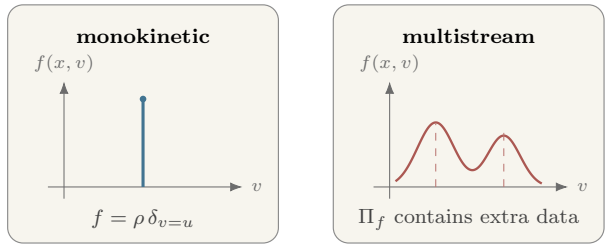
\begin{figure}[t]
\centering
\resizebox{0.98\columnwidth}{!}{%
\begin{tikzpicture}[>=Latex, font=\small]
\begin{scope}
\node[noticebox, minimum width=3.70cm, minimum height=3.25cm] at (0,0) {};
\node[font=\bfseries\footnotesize] at (0,1.20) {monokinetic};
\draw[->, kinink!70] (-1.30,-0.86) to (1.36,-0.86) node[right, font=\scriptsize] {$v$};
\draw[->, kinink!70] (-1.08,-1.02) to (-1.08,0.58) node[above, font=\scriptsize] {$f(x,v)$};
\draw[kinblue!85, line width=1.2pt] (0,-0.86) to (0,0.34);
\fill[kinblue!85] (0,0.34) circle (0.045);
\node[noticelabel] at (0,-1.30) {$f=\rho\,\delta_{v=u}$};
\end{scope}
\begin{scope}[xshift=4.45cm]
\node[noticebox, minimum width=3.70cm, minimum height=3.25cm] at (0,0) {};
\node[font=\bfseries\footnotesize] at (0,1.20) {multistream};
\draw[->, kinink!70] (-1.30,-0.86) to (1.36,-0.86) node[right, font=\scriptsize] {$v$};
\draw[->, kinink!70] (-1.08,-1.02) to (-1.08,0.58) node[above, font=\scriptsize] {$f(x,v)$};
\draw[kinred!85, line width=0.9pt, domain=-1.0:1.0, samples=120]
  plot (\x,{0.88*exp(-8*(\x+0.45)*(\x+0.45))
       +0.70*exp(-10*(\x-0.48)*(\x-0.48))-0.86});
\draw[dashed, kinred!70] (-0.45,-0.86) to (-0.45,0.00);
\draw[dashed, kinred!70] (0.48,-0.86) to (0.48,-0.10);
\node[noticelabel] at (0,-1.30) {$\Pi_f$ contains extra data};
\end{scope}
\end{tikzpicture}}
\caption{The quasineutral constraint fixes only the density, not the whole
velocity distribution.  If the distribution is concentrated at one velocity,
the moment equations close and one obtains the incompressible Euler equations.
If several velocities coexist at the same spatial point, the second moment
\(\Pi_f\) contains additional information, and the limiting equation remains
genuinely kinetic.}
\label{fig:monokinetic-multistream}
\end{figure}

This point will be important in what follows.  The constraint \(\rho_f=1\)
fixes only the integral of \(f\) in velocity; it does not determine how the
particles are distributed among the possible velocities at each point in space.
This remaining kinetic freedom is one of the sources of difficulty in the
quasineutral limit.

The rigorous analysis of this limit was initiated in the pioneering works of
Brenier and Grenier \cite{Brenier1989, Brenier2000,Grenier}.  Their results showed that the
formal passage from \((VP)_\eps\) to \((KInE)\) can be justified under suitable
stability or regularity assumptions; later work treated Penrose-stable regimes
\cite{Penrose,HanKwanRousset} and also revealed the oscillations and
instabilities that make the problem delicate.  In the next section we discuss
one way of measuring stability for such singular limits, based on distances
between phase-space distributions that are adapted to the kinetic flow.

\section{Measuring stability in phase space}

The formal limit described above suggests a clear equation, but it does not by
itself give a convergence argument.  In the quasineutral limit, the elliptic
relation in the original system degenerates into the constraint \(\rho_f=1\);
the electric field is selected through a compatibility condition, and
perturbations of the distribution are transported by the kinetic flow.  

Before proving convergence,
one therefore has to choose a notion of distance that is sensitive to the full
phase-space distribution and to the way it moves.

Strong norms are often too rigid for singular kinetic limits, because they
penalize small displacements and fast oscillations in a way that does not
reflect the transport structure of the equation.  On the other hand, one can
turn to a weaker and more geometric notion of distance, introduced in optimal
transport and known as the Monge-Kantorovich, or Wasserstein, distance.  Its
guiding idea is to compare two probability distributions by the cost of
rearranging the mass of one into the other.

Let \(\mathcal P_p(\T^d\times\R^d)\) be the set of probability measures on
phase space with finite \(p\)-th moment in velocity.  For \(p\ge1\), the
Wasserstein distance is defined by
\[
\begin{aligned}
   W_p(f,g)^p
   &=
   \inf_{\pi\in\revm{\Gamma(f,g)}}
   \int_{\left(\T^d\times\R^d\right)^2}
   c_p(z,z')\,d\pi(z,z'),\\
   c_p(z,z')
   &=
   d_{\T}(x,x')^p+|v-v'|^p,
\end{aligned}
\tag{11}
\]
where \(z=(x,v)\) and \(z'=(x',v')\).  Here \(\revm{\Gamma(f,g)}\) denotes the set of
couplings between \(f\) and \(g\), that is, the probability measures on two
copies of phase space whose marginals are \(f\) and \(g\).  One may think of a
coupling as a plan for matching the mass of \(f\) with the mass of \(g\).  \rev{The
distance \(d_{\T}\) is the periodic distance on the flat torus,
\[
   d_{\T}(x,y)=\min_{m\in\mathbb Z^d}|x-y+m|.
\]}
The
quantity in (11) is the average cost of this matching, and \(W_p\) is obtained
by choosing the cheapest plan.

This notion of distance is useful in kinetic problems for reasons that are
already visible in simple examples.  Two nearby Dirac masses are close from
the transport point of view, because
\[
   W_p(\delta_z,\delta_{z'})
   =
   \left(d_{\T}(x,x')^p+|v-v'|^p\right)^{1/p},
\]
whereas their total variation distance is maximal as soon as the two points do
not coincide.  Rapid oscillations can also be small in transport distance,
when the excess and deficit of mass alternate at a short scale.  For example,
on the one-dimensional torus, with normalized Lebesgue measure, the densities
\[
   \rho_k(x)=1+\eta\sin(\revm{2\pi} kx),\qquad 0<\eta<1,
\]
remain separated from \(1\) in \(L^1\), uniformly as \(k\to\infty\), but
converge to \(1\) in \(W_1\), because the excess and deficit of mass can be
matched over distances of order \(1/k\).  \rev{Equivalently, the measures
\(\rho_k(x)\,\dd x\) converge weakly to \(\dd x\): the oscillations disappear
when tested against any fixed continuous function, even though their \(L^1\)
size does not vanish.}  This is the kind of behavior
naturally produced by singular limits with fast oscillations.  Finally,
Wasserstein distances are compatible with weak convergence.  On compact
spaces, they metrize weak convergence of probability measures; on
\(\T^d\times\R^d\), convergence in \(W_p\) corresponds to weak convergence
together with control of the relevant velocity moments.  They therefore give a
way to measure stability that is weaker than pointwise norms, but still
quantitative and sensitive to the geometry of phase space.

This transport point of view fits naturally with Vlasov equations.  Solutions
are carried by characteristics, and a coupling between two initial
distributions can be transported by the two corresponding characteristic
flows.  In this way one obtains a coupling between the two solutions at later
times, and stability becomes a question about how fast the transportation cost
can grow along pairs of characteristics.

The use of Wasserstein distances in kinetic theory goes back to the work of
Dobrushin in the late 1970s \cite{Dobrushin}.  In the case of smooth
interaction potentials, he proved existence and uniqueness for Vlasov
equations by using precisely this transport point of view.  The key idea is to compare two solutions by matching their initial data, move
the matching by the two characteristic flows, and estimate how much the cost of
this matching can grow.

Let us write this in a simplified form.  Consider two solutions of
\[
\begin{aligned}
   \partial_t f_i+v\cdot\nabla_x f_i
   +F[f_i]\cdot\nabla_v f_i&=0,
   \qquad i=1,2,\\
   F[f_i]&=K\ast\rho_{f_i}.
\end{aligned}
\]
If the force kernel \(K\) is Lipschitz, for instance
\(\|DK\|_{L^\infty}<\infty\), then Dobrushin's argument gives
\[
   W_1(f_1(t),f_2(t))
   \le e^{Ct} W_1(f_1(0),f_2(0)).
\tag{12}
\]
Indeed, one starts from an optimal coupling of the two initial data and pushes
it forward by the two characteristic flows.  This produces an admissible
coupling of \(f_1(t)\) and \(f_2(t)\).  Along paired characteristics, the
Lipschitz regularity of the force controls the growth of the distance between
the two trajectories; minimizing over the initial coupling gives (12).

The \((VP)\) equation is more delicate because the Coulomb force is
singular, and the Lipschitz estimate used above is no longer available.
Loeper's theorem provides a fundamental replacement under bounded density
assumptions \cite{Loeper}.  In rough terms, if
\[
   A(t)=\|\rho_{f_1}(t)\|_{L^\infty\revm{(\T^d)}}
       +\|\rho_{f_2}(t)\|_{L^\infty\revm{(\T^d)}},
\]
then one obtains a stability inequality with a logarithmic correction:
\[
\begin{aligned}
   \frac{d}{dt}W_2(f_1(t),f_2(t))^2
   &\lesssim
   A(t)W_2(f_1(t),f_2(t))^2\\
   &\quad\times
   \left(1+\left|\log W_2(f_1(t),f_2(t))^2\right|\right).
\end{aligned}
\tag{13}
\]
This is an Osgood-type estimate, reflecting the log-Lipschitz regularity of
the electric field generated by a bounded density.

Loeper's estimate became a cornerstone for stability questions in
\((VP)\) theory, and it is one of the basic tools behind several
passages from kinetic to macroscopic descriptions.  In the quasineutral
problem, however, the estimate is strained by the singular scaling.  The
logarithmic correction already implies that, if the initial error is of size
\(\theta\), the useful time scale behaves roughly like
\[
   t\lesssim \log|\log\theta|.
\tag{14}
\]
At the same time, the elliptic estimate for the electric field carries the
factor \(\eps^{-2}\) coming from the scaled Poisson equation
\(-\eps^2\Delta U_\eps=\rho_{f_\eps}-1\).  Combining the Osgood loss with this
singular elliptic factor leads to very stringent smallness assumptions, such as
double-exponential closeness in \(\eps\) in higher dimensions.

The estimates above show the power of Wasserstein stability, and at the same
time explain why the classical form of the argument may become too costly in
singular limits.  In the quasineutral problem, the logarithmic loss and the
singular elliptic factor amplify small errors, leading to very strong
preparation assumptions on the initial data.  Sharper results require a more
accurate way of measuring the stability of the plasma, one that reflects more
closely how perturbations are transported in phase space.  This will be the
main topic of the next section.

\section{The geometry of kinetic transport}

To see what a more accurate distance should measure, it is useful to begin
with the simplest kinetic equation, free transport,
\[
   \partial_t f+v\cdot\nabla_x f=0.
\tag{15}
\]
Its characteristic flow is
\[
   S_t(x,v)=(x+tv,v).
\tag{16}
\]
Even in this elementary case, the geometry of phase space is tied to the
kinetic motion.  If two particles start from \((x,v)\) and \((y,w)\), their
usual phase-space distance at time \(t\), working for simplicity in \(\R^d\)
or with compatible lifts of the torus variables, is
\[
   |x-y+t(v-w)|+|v-w|.
\tag{17}
\]
Thus a small discrepancy in velocity is converted by the flow into a
discrepancy in position. A useful
stability distance should therefore distinguish the part of the displacement
which is already explained by the reference motion from the part which still
has to be estimated.

\begin{figure}[t]
\centering
\resizebox{0.98\columnwidth}{!}{%
\begin{tikzpicture}[>=Latex, font=\small]
\begin{scope}
\node[noticebox, minimum width=3.68cm, minimum height=3.55cm] at (0,0) {};
\node[font=\bfseries\footnotesize] at (0,1.22) {fixed coordinates};
\draw[->, kinink!70] (-1.38,-0.86) to (1.36,-0.86);
\node[font=\scriptsize, text=kinink!70] at (1.22,-0.72) {$x$};
\draw[->, kinink!70] (-1.12,-1.02) to (-1.12,0.62) node[above, font=\scriptsize] {$v$};
\begin{scope}[rotate=28]
\draw[fill=kinblue!13, pattern=north east lines, pattern color=kinblue!35,
  draw=kinblue!85, line width=0.75pt] (0,-0.14) ellipse (0.78 and 0.24);
\end{scope}
\node[noticelabel, text width=2.7cm] at (0,-1.34) {free transport shears the blob};
\end{scope}
\begin{scope}[xshift=4.35cm]
\node[noticebox, minimum width=3.68cm, minimum height=3.55cm] at (0,0) {};
\node[font=\bfseries\footnotesize] at (0,1.22) {adapted coordinates};
\draw[->, kinink!70] (-1.38,-0.86) to (1.36,-0.86);
\node[font=\scriptsize, text=kinink!70] at (1.10,-0.72) {$x-tv$};
\draw[->, kinink!70] (-1.12,-1.02) to (-1.12,0.62) node[above, font=\scriptsize] {$v$};
\draw[fill=kinblue!13, pattern=north east lines, pattern color=kinblue!35,
  draw=kinblue!85, line width=0.75pt] (0,-0.14) ellipse (0.62 and 0.32);
\node[noticelabel, text width=2.8cm] at (0,-1.34) {pull back by \(S_{-t}\)};
\end{scope}
\end{tikzpicture}}
\caption{Under free transport, a velocity discrepancy becomes a position
discrepancy.  Pulling back by \(S_{-t}\) removes this deterministic shear and
measures the remaining error in adapted coordinates.}
\label{fig:free-shear}
\end{figure}
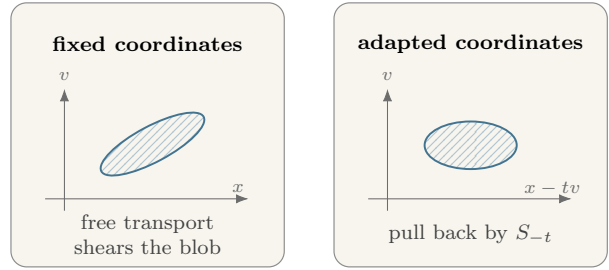

The formula above suggests a first way to calibrate the transport cost.  Instead
of comparing two distributions in the fixed coordinates \((x,v)\), one compares
them in coordinates in which the chosen reference motion has been straightened.
For \(p\ge1\), define the \rev{free-transport-adapted} distance by
\[
\begin{aligned}
   W^{\rm free}_{p,t}(f,g)
   &:=
   W_p(S_{-t}\# f,S_{-t}\# g),\\
   S_{-t}(x,v)&=(x-tv,v).
\end{aligned}
\tag{18}
\]
If \(f\) and \(g\) are both solutions of the free transport equation, then
\(W^{\rm free}_{p,t}(f(t),g(t))\) is constant in time.  In this example, the
cost has been chosen so that the reference dynamics is already built into the
distance; the stability estimate can then focus on the discrepancy which is not
explained by that dynamics.

For two point masses observed at time \(t\), this amounts to using the
quadratic cost
\[
   |(x-y)-t(v-w)|^2+|v-w|^2.
\tag{19}
\]
The first term is the residual position mismatch after subtracting the
displacement predicted by the velocity mismatch over the time interval
\([0,t]\).  This is a first example of an adapted cost in phase space; it illustrates the
general principle that, in a kinetic estimate, the cost should isolate the part
of the phase-space error that actually has to be controlled.

For a Vlasov equation, the same adapted variable is no longer conserved,
because particles are accelerated by the self-consistent field.  It still
separates the kinematic part of the motion from the part created by the force.
Let
\[
   \dot X_i=V_i,
   \qquad
   \dot V_i=E_i(t,X_i),
   \qquad i=1,2,
\]
be two characteristic systems, and set
\[
   R(t)=X_1(t)-X_2(t)-t\bigl(V_1(t)-V_2(t)\bigr).
\]
Then
\[
   \frac{d}{dt}R(t)
   =
   -t\bigl(E_1(t,X_1)-E_2(t,X_2)\bigr),
\tag{20}
\]
whereas
\[
   \frac{d}{dt}\bigl(V_1(t)-V_2(t)\bigr)
   =
   E_1(t,X_1)-E_2(t,X_2).
\tag{21}
\]
The free-transport contribution has disappeared from the evolution of
\(R(t)\).  What remains is the difference between the two fields, which is the
quantity that must be estimated in a stability argument.

This observation already gives a refinement of Dobrushin's estimate for smooth
interaction kernels.  Assume that \(K\) is \(C^1\), set
\(B=\|DK\|_{L^\infty}\), and consider the Vlasov equation with force
\(F[f]=K\ast\rho_f\).

\begin{theorem}[Kinetic Dobrushin estimate, simplified form]
Let \(f_1\) and \(f_2\) be two solutions of
\[
   \partial_t f+v\cdot\nabla_x f+F[f]\cdot\nabla_v f=0,
   \qquad
   F[f]=K\ast\rho_f,
\]
with \(K\) Lipschitz.  Then
\[
   W_1(f_1(t),f_2(t))
   \le C_B(t) W_1(f_1(0),f_2(0)),
\tag{22}
\]
where
\[
   C_B(t)=
   \min\left\{
   (1+t)e^{\frac23 B((1+t)^3-1)},
   e^{(1+2B)t}
   \right\}.
\]
\end{theorem}

\begin{ideaproof}
Start from an optimal \(W_1\)-coupling of the two initial data and transport it
by the two characteristic flows.  Instead of differentiating the usual
phase-space cost, consider
\[
   Q(t)=\int
   \Bigl(
   |X_1-tV_1-(X_2-tV_2)|+|V_1-V_2|
   \Bigr)\,d\pi_0.
\]
At time \(0\), this is the usual transportation cost, so
\(Q(0)=W_1(f_1(0),f_2(0))\).  At time \(t\), the usual distance is controlled
by \(Q(t)\) through
\[
   W_1(f_1(t),f_2(t))\le (1+t)Q(t).
\]
Differentiating \(Q(t)\) uses the cancellation in (20).  The remaining terms
involve the difference between the two forces, and the Lipschitz bound on
\(K\) gives
\[
   Q'(t)\le 2B(1+t)^2 Q(t).
\]
After integration, this yields the first term in the definition of \(C_B(t)\).
The second term is the classical Dobrushin bound, and one keeps the better of
the two estimates.
\end{ideaproof}

When \(K=0\), the quantity \(Q(t)\) is constant.  The usual phase-space
distance may grow, but in that case the growth comes only from free transport.
The adapted cost records this elementary fact and gives the optimal behavior
for the free equation.

This example is one instance of a more flexible use of Wasserstein distances
in kinetic problems.  The cost may follow a reference flow, weight position
and velocity differently, contain cross terms, or depend nonlinearly on the
size of the displacement.  There is no universal best choice; each problem
requires its own adaptation.  The common feature is that the distance is chosen
so as to reflect the asymmetry between the variables and to identify the part
of the perturbation which is relevant for the estimate.  This is the point of
view behind \rev{kinetic Wasserstein distances, namely transport distances
whose cost is adapted to the geometry of the kinetic flow}.

The same point of view has also been useful beyond the quasineutral questions
discussed here, for instance in the construction of quantum transport
pseudometrics comparing Hartree dynamics and \((VP)\)
\cite{IacobelliLafleche}, and in stability estimates for Vlasov equations in
the presence of a strong magnetic field \cite{RegeMag}.

\section{Kinetic Wasserstein stability}

The previous section described a general principle: in a kinetic estimate, the
transport cost can be chosen so as to reflect the structure of the equation.
We now return to the quasineutral problem and use this principle for the scaled
\((VP)_\eps\) system.  The difficulty is that two singular effects occur at
the same time.  The electric field is controlled through the Poisson equation,
which contains the small factor \(\eps^2\), and the field generated by a
bounded density has only logarithmic regularity, as in Loeper's estimate.
Taken together, these two effects make the classical Wasserstein estimate too
expensive for sharp quasineutral stability.

The kinetic Wasserstein distances introduced in \cite{Iacobelli2022} were
designed to use the flexibility of the transport cost in this setting.  For
\((VP)_\eps\), one keeps a quadratic phase-space cost, but allows the relative
weight of spatial and velocity displacements to depend on the quasineutral
scale and on the size of the distance itself.

Let \(f_1\) and \(f_2\) be two solutions of the scaled system \((VP)_\eps\),
and let \((X_i,V_i)\) be their characteristic flows.  Starting from an optimal
coupling \(\pi_0\) of the initial data, one considers a quantity of the form
\[
\begin{aligned}
   Q(t)
   &=
   \int
   \Bigl[
   \lambda(t)|\revm{\delta X}(t)|^2
   +|\revm{\delta V}(t)|^2
   \Bigr]\,d\pi_0,\\
   \revm{\delta X}(t)&=X_1(t)-X_2(t),\\
   \revm{\delta V}(t)&=V_1(t)-V_2(t).
\end{aligned}
\tag{23}
\]
The weight \(\lambda(t)\) is chosen so as to compensate for the singular field
estimate.  At the level of scales, the relevant choice is
\[
   \lambda(t)\sim \eps^{-2}|\log Q(t)|.
\tag{24}
\]
Thus spatial discrepancies are measured more strongly when the Poisson
equation is singular, and the logarithmic factor is tied to the modulus of
continuity of the electric field.  In the actual proof, \(Q\) is defined
implicitly in the regime where it is small; for the present discussion, the
important point is that the cost is not fixed in advance, but is adapted to the
stability estimate.

\begin{figure}[t]
\centering
\resizebox{0.98\columnwidth}{!}{%
\begin{tikzpicture}[>=Latex, font=\small]
\begin{scope}
\node[noticebox, minimum width=3.95cm, minimum height=3.34cm] at (0,0) {};
\node[font=\bfseries\footnotesize] at (0,1.20) {standard cost};
\draw[->, kinink!70] (-1.48,-0.90) to (1.48,-0.90);
\node[font=\scriptsize, text=kinink!70] at (1.23,-0.74) {$\revm{\delta X}$};
\draw[->, kinink!70] (-1.16,-1.08) to (-1.16,0.70) node[above, font=\scriptsize] {$\revm{\delta V}$};
\draw[kinink!65, line width=0.8pt] (0,-0.18) circle (0.62);
\draw[->, kinred!85, line width=0.9pt] (0,-0.18) to (0.66,0.27);
\fill[kinred!85] (0.66,0.27) circle (0.045);
\node[noticelabel] at (0,-1.38) {$|\revm{\delta X}|^2+|\revm{\delta V}|^2$};
\end{scope}
\begin{scope}[xshift=4.62cm]
\node[noticebox, minimum width=3.95cm, minimum height=3.34cm] at (0,0) {};
\node[font=\bfseries\footnotesize] at (0,1.20) {kinetic cost};
\draw[->, kinink!70] (-1.48,-0.90) to (1.48,-0.90);
\node[font=\scriptsize, text=kinink!70] at (1.23,-0.74) {$\revm{\delta X}$};
\draw[->, kinink!70] (-1.16,-1.08) to (-1.16,0.70) node[above, font=\scriptsize] {$\revm{\delta V}$};
\draw[kinblue!85, line width=0.95pt] (0,-0.18) ellipse (0.28 and 0.78);
\draw[->, kinred!85, line width=0.9pt] (0,-0.18) to (0.66,0.27);
\fill[kinred!85] (0.66,0.27) circle (0.045);
\node[noticelabel] at (0,-1.38) {$\lambda(t)|\revm{\delta X}|^2+|\revm{\delta V}|^2$};
\node[noticelabel, text=kinblue!90] at (1.05,0.66) {$\lambda(t)>1$};
\end{scope}
\end{tikzpicture}}
\caption{The adapted cost weights spatial discrepancies more strongly in the
quasineutral regime.  \rev{The red segment represents one paired phase-space
displacement, measured with two different costs.}  The weight depends both on
the singular scale \(\eps\) and on the logarithmic modulus appearing in
Loeper's estimate.}
\label{fig:anisotropic-cost}
\end{figure}
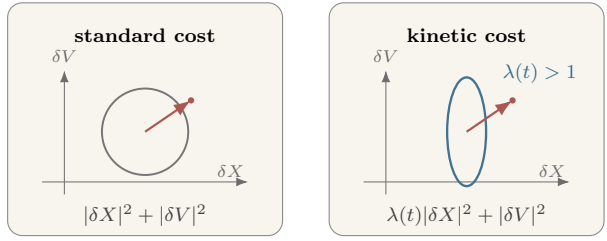

The resulting estimate can be stated schematically as follows.  The precise
constants and smallness conditions are not the main point here; what matters is
the differential inequality satisfied by the adapted quantity.

\begin{theorem}[Kinetic Loeper estimate, simplified form]
Let \(f_1\) and \(f_2\) be two solutions of \((VP)_\eps\) with bounded
densities, and set
\[
   A(t)=\|\rho_{f_1}(t)\|_{L^\infty\revm{(\T^d)}}
       +\|\rho_{f_2}(t)\|_{L^\infty\revm{(\T^d)}}.
\]
Assume that \(A\in L^1(0,T)\).  There is an adapted transport quantity \(Q(t)\)
of the form (23)-(24), well defined as long as it remains sufficiently small,
such that
\[
   \dot Q(t)
   \lesssim
   \frac{A(t)}{\eps}\,Q(t)\sqrt{|\log Q(t)|}.
\tag{25}
\]
Consequently,
\[
   Q(t)\le
   \exp\left[
   -\left(
   \sqrt{|\log Q(0)|}
   -\frac{C}{\eps}\int_0^t A(s)\,ds
   \right)_+^2
   \right].
\tag{26}
\]
\end{theorem}

This should be compared with the classical Osgood estimate discussed in the
previous section.  If the initial distance is of size \(\theta\), the
classical argument gives a useful stability time of order
\(\log|\log\theta|\).  Estimate (26) gives instead a time scale of order
\(|\log\theta|^{1/2}\).  In a quasineutral problem, where the stability loss
contains powers of \(1/\eps\), this square-root improvement is exactly what
turns double-exponential preparation thresholds into single-exponential ones.

\medskip
\noindent\textit{Deriving (26) from (25).}
The estimate (25) is useful because it becomes linear after the change of
variable
\[
   R(t)=\sqrt{|\log Q(t)|}.
\]
Indeed, as long as \(0<Q(t)<1\), we have \(Q=e^{-R^2}\), and therefore
\[
   \dot Q=-2R\dot R\,Q.
\]
Substituting this identity into (25), and dividing by \(QR\), gives
\[
   -\dot R(t)\lesssim \frac{A(t)}{\eps}.
\]
Thus
\[
   R(t)\ge R(0)-\frac{C}{\eps}\int_0^t A(s)\,ds,
\]
which is (26).  The point is that the adapted estimate gives direct control of
\(\sqrt{|\log Q|}\).  In contrast, the classical Osgood estimate controls only
\(\log|\log Q|\), and therefore requires a much smaller initial error in order
to survive the same quasineutral loss.
\medskip

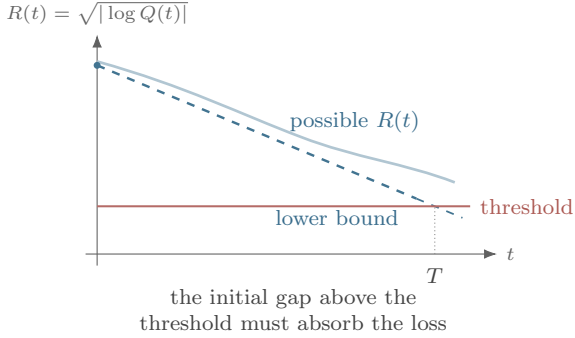
\begin{figure}[t]
\centering
\resizebox{0.96\columnwidth}{!}{%
\begin{tikzpicture}[>=Latex, font=\small]
\draw[->, kinink!75] (-0.15,0) to (5.2,0) node[right, font=\scriptsize] {$t$};
\draw[->, kinink!75] (0,-0.15) to (0,2.85)
  node[above, font=\scriptsize] {$R(t)=\sqrt{|\log Q(t)|}$};
\draw[kinblue!35, line width=1.0pt]
  (0,2.50) .. controls (1.15,2.18) and (1.85,1.82) .. (2.65,1.52)
           .. controls (3.35,1.25) and (4.00,1.18) .. (4.65,0.93);
\draw[kinblue!85, dashed, line width=0.9pt] (0,2.45) to (4.15,0.72);
\draw[kinblue!85, dashed, line width=0.65pt] (4.15,0.72) to (4.75,0.47);
\draw[kinred!75, line width=0.75pt] (0,0.62) to (4.85,0.62)
  node[right, noticelabel, text=kinred!85] {threshold};
\draw[densely dotted, kinink!55] (4.39,0) to (4.39,0.62);
\node[noticelabel] at (4.39,-0.28) {$T$};
\fill[kinblue!85] (0,2.45) circle (0.045);
\node[noticelabel, text=kinblue!90] at (3.35,1.72) {possible \(R(t)\)};
\node[noticelabel, text=kinblue!90] at (3.12,0.48) {lower bound};
\node[noticelabel, text width=4.3cm, align=center] at (2.55,-0.72)
  {the initial gap above the threshold must absorb the loss};
\end{tikzpicture}}
\caption{The kinetic estimate gives an explicit lower bound for
\(R=\sqrt{|\log Q|}\).  The dashed line represents the quantity
\(R(0)-(C/\eps)\int_0^t A(s)\,ds\) obtained from (26).  The horizontal line is a
smallness threshold: since \(R=\sqrt{|\log Q|}\), staying above this line means
that \(Q\) remains sufficiently small.  To keep the adapted distance small up
to time \(T\), the initial gap above the threshold must dominate the loss
accumulated along the estimate, namely \((C/\eps)\int_0^T A(s)\,ds\).}
\label{fig:sqrt-log-gain}
\end{figure}

We now explain how this estimate is used in the quasineutral limit.  Suppose
that \(g_{0,\eps}\) is a regular family of initial data for which the
\rev{quasineutral limit is already known after removing the fast plasma
oscillations present in this regime.  This standard removal of oscillations is
often called filtering.}  Let \(g_\eps(t)\) be the corresponding solutions.  We
then perturb the initial data and write
\[
   f_{0,\eps}=g_{0,\eps}+h_{0,\eps},
\]
where the perturbation \(h_{0,\eps}\) may be rough.  The question is how small
this perturbation has to be in Wasserstein distance in order for the solution
\(f_\eps\) to have the same filtered limit as \(g_\eps\).

This rough-perturbation strategy was developed, in joint work with Han-Kwan,
first \rev{in one dimension, for both the fixed-background \((VP)\) model and
VPME,} and then for higher-dimensional \((VP)\) regimes
\cite{HanKwanIacobelli1,HanKwanIacobelli2}.  The statement below is a
schematic form of this line of results; it is meant to isolate the stability
mechanism, rather than to summarize all quasineutral limit theorems.

The quantity that has to be paid in the stability estimate is
\[
   M_\eps(T)
   :=
   \frac1\eps\int_0^T
   \bigl(
   \|\rho_{f_\eps}(t)\|_{L^\infty\revm{(\T^d)}}
   +\|\rho_{g_\eps}(t)\|_{L^\infty\revm{(\T^d)}}
   \bigr)\,dt.
\tag{27}
\]
It contains the singular factor \(1/\eps\) and the density bounds entering the
Loeper estimate.  In applications one often obtains a bound of the form
\(M_\eps(T)\lesssim \eps^{-\zeta}\), where the exponent \(\zeta\) includes the
losses coming from the quasineutral scaling and from the available a priori
estimates.

\begin{theorem}[Rough perturbations near a stable quasineutral limit]
Assume that the reference family \(g_\eps\) has a known filtered quasineutral
limit on \([0,T]\), and that the comparison with a second family \(f_\eps\)
satisfies \(M_\eps(T)<\infty\).  If, for some sufficiently large constant
\(\kappa\),
\[
   W_2(f_{0,\eps},g_{0,\eps})
   \le \exp\bigl(-\kappa M_\eps(T)^2\bigr),
\tag{28}
\]
then the adapted distance between \(f_\eps(t)\) and \(g_\eps(t)\) stays small
on \([0,T]\).  Consequently \(f_\eps\) has the same filtered quasineutral limit
as \(g_\eps\).
\end{theorem}

If \(M_\eps(T)\lesssim \eps^{-\zeta}\), condition (28) becomes the
single-exponential smallness assumption
\[
   W_2(f_{0,\eps},g_{0,\eps})
   \le \exp(-\kappa\eps^{-2\zeta}).
\]
Indeed, this gives
\[
   \sqrt{\left|\log W_2(f_{0,\eps},g_{0,\eps})^2\right|}
   \gtrsim \sqrt\kappa\, M_\eps(T),
\]
which is the logarithmic margin required by (26).  Polynomial smallness would
give only \(\sqrt{|\log \eps|}\) on the left-hand side, which is too small when
the stability loss is a negative power of \(\eps\).

This positive result should be read together with \rev{instability results.  In
particular, Han-Kwan and Hauray \cite{HanKwanHauray} showed, already for the
one-dimensional \((VP)\) equation, that near unstable profiles polynomial-size
perturbations may destroy the quasineutral limit.  The same mechanism can be
viewed in higher dimensions by considering data depending on one spatial
variable.}  The perturbative result above has a different nature: it is a
weak-strong statement around a stable reference family.  The
regular family \(g_\eps\) carries the strong information, such as analytic or
Penrose-stable structure, density bounds, and the known filtered limit.  The
family \(f_\eps\) may be much rougher, but it inherits the same limit because
the kinetic Wasserstein estimate keeps it close to \(g_\eps\) in an adapted
distance.

In this way, the kinetic Wasserstein estimate refines the classical stability
theory of Dobrushin and Loeper in the direction needed for quasineutral
limits.  It does not change the limiting equation, but it changes the class of
perturbations for which the limiting process can be justified.

\section{VPME and the electromagnetic horizon}

There is another electrostatic model which is very natural physically and more
delicate mathematically.  It is usually called VPME, for \((VP)\) with
massless, or thermalized, electrons.  In this model the ions are described
kinetically, while the electrons are assumed to remain in thermal equilibrium.
The electrostatic potential then solves a nonlinear Poisson-Boltzmann
equation:
\[
\left\{
\begin{aligned}
\partial_t f_\eps&+v\cdot\nabla_x f_\eps
    +E_\eps\cdot\nabla_v f_\eps=0,\\
E_\eps&=-\nabla_x U_\eps,\\
\eps^2\Delta_x U_\eps&=e^{U_\eps}-\rho_{f_\eps}.
\end{aligned}
\right.
\tag{29}
\]
The term \(e^{U_\eps}\) is the Maxwell-Boltzmann electron density.  VPME
appears in the physics of ion-acoustic waves, plasma expansion into vacuum,
and related phenomena; see for instance the classical work of Gurevich and
Pitaevskii \cite{GurevichPitaevskii}.  The terminology ``massless electrons''
reflects the physical separation of time scales: electrons are much lighter
than ions and are modeled as equilibrating much faster.

The formal quasineutral limit of (29) gives
\[
   e^U=\rho_f,
   \qquad
   U=\log\rho_f.
\]
Hence the limiting kinetic equation is
\[
   \partial_t f+v\cdot\nabla_x f-\nabla_x\log\rho_f\cdot\nabla_v f=0.
\tag{30}
\]
This is the kinetic isothermal Euler equation.  If one imposes the
monokinetic ansatz \(f=\rho\,\delta_{v=u}\), then the moment equations close
and give
\[
   \partial_t\rho+\nabla_x\cdot(\rho u)=0,
   \qquad
   \partial_t(\rho u)+\nabla_x\cdot(\rho u\otimes u)=-\nabla_x\rho.
\tag{31}
\]
This is the isothermal Euler system.  In the kinetic equation (30), however,
the unknown remains the full distribution \(f(t,x,v)\).

VPME is more delicate than the fixed-background \((VP)\) model for a
concrete reason: the elliptic equation is nonlinear.  The potential
\(U_\eps\) appears inside the exponential \(e^{U_\eps}\), and \rev{this couples
the field, the thermalized electron density, and the growth of the
characteristics.  Thus the stability problem is no longer governed only by the
Coulomb field estimate; one also has to control the Poisson-Boltzmann equation
and how this control feeds back into the particle trajectories.}

\rev{For this VPME problem, joint work with Griffin-Pickering shows that the
kinetic Wasserstein viewpoint gives the correct stability scale, but it has to
be combined with estimates specific to the nonlinear Poisson-Boltzmann
coupling \cite{GPI2020,GPI2026}.}  In the most recent form of the result, the
perturbation is assumed to be single-exponentially small in Wasserstein
distance,
\[
   W_1(f_{0,\eps},g_{0,\eps})
   \le \exp(-\kappa\eps^{-\alpha}),
\tag{32}
\]
where the exponent \(\alpha\) depends on the dimension and on the available
moment bounds.  For the present discussion, the scale of the assumption is the
\rev{main feature}.  Earlier arguments in this perturbative program led, in
some regimes, to double-exponential thresholds for the fixed-background
\((VP)\) model and to quadruple-exponential thresholds for VPME.  These were
thresholds imposed by the available stability estimates, rather than universal
labels of the models themselves.  \rev{The refined VPME analysis reduces the
required closeness to a single exponential scale.}

\begin{table}[t]
\centering
\footnotesize
\begin{tabularx}{\columnwidth}{@{}>{\raggedright\arraybackslash}p{0.25\columnwidth}>{\raggedright\arraybackslash}X@{}}
\toprule
Model & Positive stability threshold near a stable reference family\\
\midrule
VP, $d=1$
& sharp single-exponential scale in Wasserstein distance
  \cite{HanKwanIacobelli1}\\
\addlinespace
VPME, $d=1$
& earlier double-exponential threshold \cite{HanKwanIacobelli1};
  refined VPME stability gives a single-exponential threshold
  \cite{GPI2026}\\
\addlinespace
VP, $d=2,3$
& earlier double-exponential threshold \cite{HanKwanIacobelli2};
  kinetic Wasserstein stability gives a single-exponential threshold,
  \rev{with exponents determined by the available a priori bounds}
  \cite{Iacobelli2022}\\
\addlinespace
VPME, $d=2,3$
& earlier quadruple-exponential threshold \cite{GPI2020};
  refined VPME stability gives a single-exponential threshold, with
  dimension-dependent exponents \cite{GPI2026}\\
\bottomrule
\end{tabularx}
\caption{Selected positive thresholds in the rough-perturbation program for
quasineutral limits.  The table summarizes this particular line of work by
Han-Kwan-Iacobelli and Griffin-Pickering-Iacobelli; broader quasineutral
limit results are outside its scope.  Near unstable profiles, polynomial-size
perturbations may fail: for VP this is the instability result of
Han-Kwan-Hauray \cite{HanKwanHauray}, while for VPME the corresponding
instability mechanism is treated in the refined VPME work of
Griffin-Pickering-Iacobelli \cite{GPI2026}.}
\label{tab:thresholds}
\end{table}

The discussion so far has stayed mostly within electrostatic models.  The next
natural challenge is electromagnetic.  The relativistic Vlasov-Maxwell system
describes charged particles coupled to the full Maxwell equations.  In the
quasineutral regime the electric field has both a longitudinal component,
related to charge separation, and a solenoidal component coupled to the
magnetic field.  Unlike in the electrostatic case, the field has its own wave
dynamics.  In particular, transverse electromagnetic oscillations have to be
separated from the longitudinal plasma oscillations already present in
\((VP)\).

Recent work with Gagnebin, Rege, and Rossi proves a quasineutral limit for the
relativistic Vlasov-Maxwell system in analytic regularity, after introducing
appropriate electromagnetic correctors and \rev{filtering the electromagnetic
oscillations}
\cite{GagnebinIacobelliRegeRossi}.  The limiting model is a kinetic version of
electron magnetohydrodynamics, in a scaling where magnetic dynamics survive in
the limit.  This result is of a different nature from the rough Wasserstein
stability theorems discussed above, but it points to the same broad
mathematical issue: \rev{singular plasma limits require estimates that respect
both the transport structure of the particles and the wave structure of the
fields.}

The examples above point toward the electromagnetic problem.  In
Vlasov-Maxwell, the quasineutral limit has to be combined with the dynamics of
the electromagnetic field itself; longitudinal plasma oscillations, transverse
waves, and magnetic effects interact in the same singular regime.  The analytic
result mentioned above shows that this passage can be justified when suitable
correctors and filtrations are introduced, but it also indicates that the
electromagnetic setting will require new techniques beyond the electrostatic
ones.

Kinetic Wasserstein distances remain a useful viewpoint in this broader
program.  They show that, in singular kinetic limits, stability is often tied
to the way one measures the error in phase space.  For \((VP)\) this leads to
distances adapted to the kinetic transport and to the singular elliptic field.
For VPME, the same viewpoint has to be combined with estimates
for the nonlinear Poisson-Boltzmann coupling.  For Vlasov-Maxwell, one
expects that the relevant notions of stability will also have to interact with
the propagation of electromagnetic waves.  Understanding these mechanisms is
one of the main challenges in the analysis of quasineutral plasma models.

\end{document}